\documentclass[11pt]{article}
\usepackage{amssymb,latexsym}
\usepackage{epsfig}
\usepackage{eufrak}
\usepackage{amsmath}
\usepackage{mathrsfs}
\usepackage{color}

\newtheorem{theorem}{Theorem}
\newtheorem{lemma}{Lemma}

\newcommand{\be}{\begin{equation}}
\newcommand{\ee}{\end{equation}}
\newcommand{\bee}{\begin{eqnarray*}}
\newcommand{\eee}{\end{eqnarray*}}
\newcommand{\bel}{\begin{eqnarray}}
\newcommand{\eel}{\end{eqnarray}}
\newcommand{\bec}{\begin{cases}}
\newcommand{\eec}{\end{cases}}
\newcommand{\bem}{\begin{bmatrix}}
\newcommand{\eem}{\end{bmatrix}}

\newcommand{\la}{\label}
\newcommand{\li}{\left}
\newcommand{\ri}{\right}

\newcommand{\ovl}{\overline}

\newcommand{\Ga}{\Gamma}

\newcommand{\se}{\theta}

\newcommand{\al}{\alpha}

\newcommand{\ka}{\kappa}

\newcommand{\Om}{\Omega}

\newcommand{\f}{\frac}

\newcommand{\cd}{\cdots}

\newcommand{\qu}{\quad}
\newcommand{\qqu}{\qquad}

\newcommand{\mscr}{\mathscr}
\newcommand{\mcal}{\mathcal}

\newcommand{\bb}{\mathbb}

\newcommand{\bs}{\boldsymbol}

\newcommand{\tx}{\text}

\newcommand{\pa}{\partial}

\newcommand{\bed}{\begin{description}}
\newcommand{\eed}{\end{description}}
\newcommand{\bei}{\begin{itemize}}
\newcommand{\eei}{\end{itemize}}
\newcommand{\ben}{\begin{enumerate}}
\newcommand{\een}{\end{enumerate}}
\newcommand{\bib}{\bibitem}
\newcommand{\beL}{\begin{lemma}}
\newcommand{\eeL}{\end{lemma}}
\newcommand{\beT}{\begin{theorem}}
\newcommand{\eeT}{\end{theorem}}
\newcommand{\sect}{\section}

\newcommand{\bpf}{\begin{pf}}
\newcommand{\epf}{\end{pf}}
\newcommand{\bsk}{\bigskip}

\newcommand{\pfbox}{\hfill\mbox{$\Box$}}

\newenvironment{pf}{\paragraph*{Proof{\rm.}}}{\pfbox\bigskip}

\begin{document}

\title{{\bf Concentration Inequalities for Bounded Random Vectors}
\thanks{The author is currently working with the Department of Electrical Engineering and Computer Sciences at Louisiana
State University,  Baton Rouge, LA 70803, USA, and the Department of Electrical Engineering at Southern University and A\&M College, Baton
Rouge, LA 70813, USA; Email: chenxinjia@gmail.com}}

\author{Xinjia Chen}

\date{August 30,  2013}

\maketitle

\begin{abstract}

We derive simple concentration inequalities for bounded random vectors, which generalize Hoeffding's inequalities for bounded scalar random
variables.  As applications, we apply the general results to multinomial and Dirichlet distributions to obtain multivariate concentration
inequalities.

\end{abstract}

\section{Introduction}

Concentration phenomena is one of the most important subjects in probability theory. Formally, let $\bs{\mcal{X}}_1, \cd, \bs{\mcal{X}}_n$ be
independent random vectors of the same dimension defined in a probability space $(\Om, \Pr, \mscr{F})$. Define \be \la{defave89}
\ovl{\bs{\mcal{X}}}_n = \f{ \sum_{i=1}^n \bs{\mcal{X}}_i }{n} \ee and mathematical expectation \be \la{mudef888}
 \bs{\mu} = \bb{E} [
\ovl{\bs{\mcal{X}}}_n ]. \ee It is a fundamental problem to investigate how $\ovl{\bs{\mcal{X}}}_n$ deviate from the mean vector $\bs{\mu}$. In
many situations,
 since the exact distribution of $\ovl{\bs{\mcal{X}}}_n$ is not readily tractable, it is desirable to obtain upper bounds for probabilities
 \be
 \la{pro889}
\Pr \{ \ovl{\bs{\mcal{X}}}_n \leq \bs{z} \}   \qu \tx{and} \qu \Pr \{ \ovl{\bs{\mcal{X}}}_n \geq \bs{z} \},
 \ee
where $\bs{z}$ is a deterministic vector of the same dimension as $\ovl{\bs{\mcal{X}}}_n$, the symbols ``$\leq$'' and ``$\geq$'' are used to
denote the partial order relations between two vectors with the following notions:

For two vectors  $\bs{x} = [x_1, \cd, x_\ka]^\top$ and $\bs{y} = [y_1, \cd, y_\ka]^\top$ of the same dimension $\ka$, where $\top$ denotes
transpose operation, we write $\bs{x} \leq \bs{y}$ if
\[
x_i \leq y_i \qu \tx{for} \; i = 1, \cd, \ka.
\]
Similarly, we write $\bs{x} \geq \bs{y}$ if
\[
x_i \geq y_i \qu \tx{for} \; i = 1, \cd, \ka.
\]
Throughout this paper, such notations will be used for denoting the relations between vectors.

 To obtain upper bounds for probabilities in
(\ref{pro889}), we need some constraints for the random vectors. Actually, in many applications, the relevant random vectors are bounded.  More
specifically, for random vectors $\bs{\mcal{X}}_i = [ X_{1, i}, X_{2, i}, \cd, X_{\ka, i} ]^\top, \; i = 1, \cd, n$, there exist real numbers
$a_{\ell, i}, \; b_{\ell, i}$  such that
\[
a_{\ell, i} \leq X_{\ell, i} \leq b_{\ell, i}
\]
for $i = 1, \cd, n$ and  $\ell = 1, \cd, \ka$.  If we define
\[
Y_{\ell, i} = \f{ X_{\ell, i} - a_{\ell, i} }{ \sum_{\ell = 1}^\ka (  b_{\ell, i} - a_{\ell, i} )},
\]
then
\[
\sum_{\ell = 1}^k Y_{\ell, i} \leq 1, \qqu i = 1, \cd, n
\]
and
\[
Y_{\ell, i} \geq 0, \qqu 1 \leq \ell \leq k, \qqu 1 \leq i \leq n.
\]
It follows that random vectors $\bs{\mcal{Y}}_i = [ Y_{1, i}, Y_{2, i}, \cd, Y_{\ka, i} ]^\top, \; i = 1, \cd, n$ are bounded in a simplex. This
demonstrates that a bounded random vector can be transformed into a random vector contained in a simplex by virtue of the translation and
scaling operations.  Motivated by this analysis, we shall  derive upper bounds for the probabilities in (\ref{pro889}) with relevant random
vectors bounded in a simplex.

The remainder of the paper is organized as follows.  In Section 2, we present our main results for multivariate concentration inequalities for
bounded random vectors.  In Section 3, we apply our general results to multinomial and Dirichlet distributions to obtain concentration
inequalities.  The proof of our main results is given in Section 4.

\section{Main Results}

Let $\bs{\mcal{X}}_i = [ X_{1, i}, X_{2, i}, \cd, X_{\ka, i} ]^\top, \; i = 1, \cd, n$ be independent random vectors  of dimension $k$ such that
\[
\sum_{\ell = 1}^k X_{\ell, i} \leq 1 \qu \tx{for} \; i = 1, \cd, n
\]
and
\[
X_{\ell, i} \geq 0 \qu \tx{for} \; 1 \leq \ell \leq k, \;  1 \leq i \leq n.
\]
Note that the elements $X_{1, i}, X_{2, i}, \cd, X_{\ka, i}$ of $\bs{\mcal{X}}_i$ are not necessarily independent.  Let $\ovl{\bs{\mcal{X}}}_n$
be the average random vector defined by (\ref{defave89}). Then, the mean $\bs{\mu}$ of $\ovl{\bs{\mcal{X}}}_n$ in accordance with
(\ref{mudef888}) can be expressed as
\[ \bs{\mu} = [\mu_1, \cd, \mu_\ka]^\top,
\]
where
\[ \mu_\ell = \f{ \sum_{i=1}^n
\bb{E} [ X_{\ell, i} ] }{n}, \qqu \ell = 1, \cd, \ka.
\]
In this setting, we have established the following results.

\beT

\la{general} Let $z_\ell, \; \ell = 0, 1, \cd, \ka$ be positive real numbers such that $\sum_{\ell = 0}^k z_\ell= 1$. Define
\[
\mu_0 = 1 - \sum_{\ell = 1}^\ka \mu_\ell, \qqu z_0 = 1 - \sum_{\ell = 1}^\ka z_\ell, \qqu  \bs{z} = [z_1, \cd, z_\ka]^\top. \] Then,  \[ \Pr \{
\ovl{\bs{\mcal{X}}}_n \leq \bs{z} \} \leq \prod_{\ell = 0}^k \li ( \f { \mu_\ell  }{z_\ell  } \ri )^{n z_\ell} \qqu \tx{provided that} \; \bs{z}
\leq \bs{\mu}.
\]
Similarly, \[ \Pr \{ \ovl{\bs{\mcal{X}}}_n \geq \bs{z} \} \leq \prod_{\ell = 0}^k \li ( \f { \mu_\ell  }{z_\ell  } \ri )^{n z_\ell} \qqu
\tx{provided that} \; \bs{z} \geq \bs{\mu}. \]

\eeT

It should be noted that Theorem \ref{general} is a multivariate generalization of Hoeffding's inequality \cite{Hoeffding}.

\section{Applications}

In this section, we shall apply the main results to multinomial distribution and Dirichlet distribution.

\subsection{Multinomial Distribution}

In probability theory, random variables $X_0, X_1, \cd, X_\ka$ are said to possess a multinomial distribution if they have a probability mass
function
 \be \la{multnomial} \Pr \{ X_i = x_i, \; i = 0, 1, \cd, \ka \}  = \f{n!}{ \prod_{i=0}^k x_i! } \prod_{i=0}^k {p_i}^{x_i},  \ee
 where $x_0, x_1, \cd, x_\ka$ are non-negative integers sum to $1$ and $p_i, \; i = 0, 1, \cd, \ka$ are positive real numbers sum to $1$.

Define $\bs{X} = [X_1, \cd,  X_\ka ]^\top$ and $\bs{\mu} = \bb{E} [ \bs{X} ]$.  Then, $\bs{\mu} = [\mu_1, \cd, \mu_\ka]^\top$, where \[ \mu_i =
\bb{E} [ X_i ] = n p_i, \qqu i = 0, 1, \cd, \ka.
\]
Let $\bs{\mcal{Y}} = [Y_1, \cd, Y_\ka]^\top$ be a random vector such that \[ \Pr \{ Y_i = y_i, \; i = 1, \cd, \ka \}  =  \prod_{i=0}^k
{p_i}^{y_i},  \]
 where $y_0, y_1, \cd, y_\ka$ are non-negative integers sum to $1$. Let $\bs{\mcal{Y}}_1, \cd \bs{\mcal{Y}}_n$ be independent random vectors
 having the same distributions as $\bs{\mcal{Y}}$.
 It is well known that random vector $\bs{X}$ has the same distribution as  $\sum_{i=1}^n \bs{\mcal{Y}}_i$. This property is referred to as the reproducibility of the multinomial
 distribution.  As applications of Theorem \ref{general} and the reproducibility of  the multinomial distribution, we have the following results.

\beT

\la{Genmulmultinomial}

Let $z_i, \; i = 0, 1, \cd, k$ be nonnegative integers such that $\sum_{i = 0}^k z_i= n$.  Let $\bs{z} = [z_1, \cd, z_\ka]^\top$. Then, \be
\la{inevipnomiala}
 \Pr \{ \bs{X} \leq \bs{z} \} \leq \prod_{i = 0}^k \li ( \f{ \mu_i } { z_i } \ri)^{z_i}  \qu \tx{provided that
$\bs{z} \leq  \bs{\mu}$}, \ee and \be \la{inevipnomialb}
 \Pr \{ \bs{X} \geq \bs{z} \} \leq \prod_{i = 0}^k \li ( \f{ \mu_i } { z_i } \ri)^{z_i} \qu \tx{provided that $\bs{z} \geq \bs{\mu}$}.  \ee

\eeT

It should be noted that Theorem \ref{Genmulmultinomial} is a multivariate generalization of the Chernoff-Hoeffding bounds for binomial
distributions \cite{Chernoff, Hoeffding}. The results in Theorem \ref{Genmulmultinomial} had been established by Chen \cite[page 17, Corollary
1]{Chen2013} with a likelihood ratio method.

\subsection{Dirichlet Distribution}

In probability theory, random variables $X_0, X_1, \cd, X_\ka$ are said to possess a Dirichlet distribution if they have a probability density
function \[ f (\bs{x}, \bs{\al})  = \bec \f{1}{ \mcal{B} (\bs{\al})} \prod_{i=0}^\ka
x_i^{\al_i - 1} & \tx{for $x_i \geq 0, \; i = 0, 1\cd, \ka$ such that  $\sum_{i=0}^\ka x_i = 1$},\\
0 & \tx{else} \eec
\]
where $\bs{x} = [x_0, x_1, \cd, x_\ka]^\top$, $\al_0, \al_1, \cd, \al_\ka$ are positive real numbers, $\bs{\al} = [ \al_0, \al_1, \cd, \al_\ka
]^\top$, and
\[
\mcal{B} (\bs{\al}) =  \f{ \prod_{i=0}^\ka \Ga (\al_i) }{ \Ga ( \sum_{i=0}^\ka \al_i)  },
\]
with $\Ga(.)$ representing the Gamma function.  The means of $X_0, X_1, \cd, X_\ka$ are
\[
\mu_i = \bb{E} [X_i] = \f{\al_i}{\sum_{\ell = 0}^\ka \al_\ell} \qqu \tx{for $i = 0, 1, \cd, \ka$}.
\]
Define vectors
\[
\bs{X} = [X_1, \cd, X_\ka]^\top, \qqu \bs{\mu} = [\mu_1, \cd, \mu_\ka]^\top.
\]
Let $\bs{\mcal{X}}_1, \bs{\mcal{X}}_2, \cd, \bs{\mcal{X}}_n$ be independent random vectors possessing the same distribution as $\bs{X}$.  Define
$\ovl{\bs{\mcal{X}}}_n = \f{ \sum_{\ell = 1}^n \bs{\mcal{X}}_i }{n}$.  As applications of Theorem \ref{general}, we have the following results.

\beT

\la{mulDrich} Let $z_\ell, \; \ell = 0, 1, \cd, \ka$ be positive real numbers such that $\sum_{\ell = 0}^k z_\ell= 1$. Define $\bs{z} = [z_1,
\cd, z_\ka]^\top$.  Then, \be \la{conthm1388} \Pr \{ \ovl{\bs{\mcal{X}}}_n  \leq \bs{z} \} \leq \li [ \prod_{i=0}^\ka  \li ( \f{\mu_i}{z_i} \ri
)^{z_i}  \ri ]^n \qqu \tx{provided that $\bs{z} \leq \bs{\mu}$}, \ee \be \la{conthm1388b} \Pr \{ \ovl{\bs{\mcal{X}}}_n \geq \bs{z} \} \leq \li [
\prod_{i=0}^\ka  \li ( \f{\mu_i}{z_i} \ri )^{z_i}  \ri ]^n \qqu \tx{provided that $\bs{z} \geq \bs{\mu}$}. \ee

\eeT

\sect{Proof of Theorem \ref{general}}

We need some preliminary results.

\beL Assume that $\sum_{\ell = 1}^k x_\ell \leq 1$ and $x_\ell \geq 0$ for $\ell = 1, \cd, k$.  Then, \be \la{ineq966}
 \prod_{\ell = 1}^k \exp
\li ( t_\ell  x_\ell \ri ) \leq 1 - \sum_{\ell = 1}^k x_\ell + \sum_{\ell = 1}^k x_\ell \exp ( t_\ell) \ee holds for arbitrary real numbers
$t_\ell, \; \ell = 1, \cd, k$.

\eeL

\bpf

We shall use a probabilistic approach.  Define $t_0 = 0$ and
\[
x_ 0 = 1 - \sum_{\ell = 1}^k x_\ell.
\]
Since $\sum_{\ell = 0}^k x_\ell = 1$ and $0 \leq x_\ell \leq 1$ for $\ell = 0, 1, \cd, k$, we can define a random variable $T$ such that
\[
\Pr \{ T = t_\ell \} = x_\ell, \qqu \ell = 0, 1, \cd, k.
\]
Then, the expectation of $e^T$ is \be \la{use889a}
 \bb{E} [  e^T ] = \sum_{\ell = 0}^k x_\ell \exp ( t_\ell).
\ee Note that \be \la{use889b} \exp ( \bb{E} [ T ] ) = \exp \li ( \sum_{\ell = 0}^k t_\ell  x_\ell \ri ) = \prod_{\ell = 0}^k \exp \li ( t_\ell
x_\ell \ri ). \ee
 Since $e^t$ is a convex function of real number $t$, it follows from Jensen's inequality that
 \be \la{use889c}
 \bb{E} [  e^T ] \geq \exp ( \bb{E} [ T ] ).
 \ee
Making use of (\ref{use889a}),  (\ref{use889b}) and (\ref{use889c}) yields, \be \la{revise88}
 \sum_{\ell = 0}^k x_\ell \exp ( t_\ell) \geq
\prod_{\ell = 0}^k \exp \li ( t_\ell x_\ell \ri ). \ee Since $t_0 = 0$ and $x_ 0 = 1 - \sum_{\ell = 1}^k x_\ell$, the inequality
(\ref{revise88}) can be written as (\ref{ineq966}).  This completes the proof of the lemma.

\epf

\beL \la{howgood}
\[
\prod_{i=1}^n \bb{E} \li [  \prod_{\ell = 1}^k \exp \li ( t_\ell  X_{\ell, i} \ri ) \ri ] \leq \li [ \mu_0 + \sum_{\ell = 1}^k \mu_\ell \exp (
t_\ell ) \ri ]^n
\]
for arbitrary real numbers $t_\ell, \; \ell = 1, \cd, k$.

\eeL

\bpf

By the independence of the random vectors $\bs{\mcal{X}}_1, \cd, \bs{\mcal{X}}_n$, we have that \be \la{indep88}
 \bb{E} \li [  \prod_{\ell = 1}^k \exp \li ( t_\ell \sum_{i=1}^n X_{\ell, i} \ri ) \ri ] = \prod_{i=1}^n \bb{E} \li [
\prod_{\ell = 1}^k \exp \li ( t_\ell  X_{\ell, i} \ri ) \ri ] \ee holds for arbitrary real numbers $t_\ell, \; \ell = 1, \cd, k$.  By Lemma 1,
we have that
\[
\prod_{\ell = 1}^k \exp \li ( t_\ell  X_{\ell, i} \ri ) \leq 1 - \sum_{\ell = 1}^k X_{\ell, i} + \sum_{\ell = 1}^k X_{\ell, i} \exp ( t_\ell),
\qqu i = 1, \cd, n
\]
holds for arbitrary real numbers $t_\ell, \; \ell = 1, \cd, k$.    Taking expectation on both sides of the above inequality yields \be
\la{prod888}
 \bb{E} \li [ \prod_{\ell = 1}^k \exp \li ( t_\ell  X_{\ell, i} \ri ) \ri ] \leq 1 - \sum_{\ell
= 1}^k \bb{E} [ X_{\ell, i} ] + \sum_{\ell = 1}^k \bb{E} [ X_{\ell, i} ] \exp ( t_\ell), \qqu i = 1, \cd, n. \ee It follows from (\ref{indep88})
and (\ref{prod888}) that \bel \bb{E} \li [  \prod_{\ell = 1}^k \exp \li ( t_\ell \sum_{i=1}^n X_{\ell, i} \ri ) \ri ] & =  & \prod_{i=1}^n
\bb{E}
\li [ \prod_{\ell = 1}^k \exp \li ( t_\ell  X_{\ell, i} \ri ) \ri ] \nonumber\\
& \leq & \prod_{i=1}^n \li [ 1 -
\sum_{\ell = 1}^k \bb{E} [ X_{\ell, i} ] + \sum_{\ell = 1}^k \bb{E} [ X_{\ell, i} ] \exp ( t_\ell) \ri ] \nonumber\\
& \leq &  \li [ 1 - \sum_{\ell = 1}^k \mu_\ell + \sum_{\ell = 1}^k \mu_\ell \exp ( t_\ell ) \ri ]^n \la{citeineq}\\
& = & \li [ \mu_0 + \sum_{\ell = 1}^k \mu_\ell \exp ( t_\ell ) \ri ]^n,  \nonumber \eel where the inequality (\ref{citeineq}) follows from the
fact that the geometric mean does not exceed the arithmetic mean.  This completes the proof of the lemma.

\epf

\beL

\la{def68689}

Define \be \la{def86636}
 \mscr{M} (t_1, \cd, t_k) =  - \sum_{\ell = 1}^k t_\ell z_\ell  +  \ln \li [ \mu_0 + \sum_{\ell = 1}^k \mu_\ell \exp (
t_\ell ) \ri ]. \ee Then, $\f{ \pa \mscr{M} (t_1, \cd, t_k)  }{\pa t_\ell} = 0$ for $\ell = 1, \cd, k$ if and only if \[ t_\ell = \ln \f{z_\ell
\mu_0}{z_0 \mu_\ell }, \qqu \ell = 1, \cd, k.
\]

\eeL

\bpf

It can be checked that the partial derivatives of $\mscr{M} (t_1, \cd, t_k)$ are given as
\[
\f{ \pa \mscr{M} (t_1, \cd, t_k)  }{\pa t_\ell} = - z_\ell + \f{\mu_\ell \exp ( t_\ell )}{ \mu_0 + \sum_{j = 1}^k \mu_j \exp ( t_j ) }, \qqu
\ell = 1, \cd, k.
\]
Letting the  partial derivatives to be $0$ yields \be \la{letzero}
 - z_\ell + \f{\mu_\ell \exp ( t_\ell )}{ \mu_0 + \sum_{j = 1}^k \mu_j \exp ( t_j ) } = 0, \qqu
\ell = 1, \cd, k. \ee For simplicity of notations, define \be \la{def8996} \se_\ell = \mu_\ell \exp ( t_\ell ),  \qqu \ell = 1, \cd, k. \ee
Then, we can write (\ref{letzero}) as \be \la{sim8996}
 - z_\ell + \f{\se_\ell}{ \mu_0 + \sum_{j = 1}^k \se_j } = 0, \qqu
\ell = 1, \cd, k. \ee From equation (\ref{sim8996}), we have \be \la{ratio899}
 \f{\se_\ell}{ \se_i } = \f{z_\ell}{z_i}, \qqu \tx{for} \; \ell,
\; i \in \{ 1, \cd, k \}. \ee Making use of (\ref{sim8996}) and (\ref{ratio899}), we have \be \la{con8869}
 - z_\ell + \f{1}{ \f{ \mu_0
}{\se_\ell} + \sum_{j = 1}^k \f{z_j}{z_\ell} } = 0, \qqu \ell = 1, \cd, k. \ee  Since $z_\ell > 0, \; \ell = 1, \cd, k$, we can write
(\ref{con8869}) as \be \la{tt88936}
 - 1 + \f{1}{ \f{z_\ell \mu_0 }{\se_\ell} + \sum_{j = 1}^k z_j } = 0, \qqu \ell = 1, \cd, k.
\ee Making use of (\ref{tt88936}) and the fact that $z_0= 1 - \sum_{j = 1}^k z_j$, we have \be \la{solution}
 \se_\ell = \f{z_\ell \mu_0 } { z_0
}, \qqu \ell = 1, \cd, k. \ee Making use of  (\ref{def8996}) and (\ref{solution}) yields $t_\ell = \ln \f{z_\ell \mu_0}{z_0 \mu_\ell }, \;  \ell
= 1, \cd, k$.  The lemma is thus proved.  \epf

\bsk

We are now in a position to prove the theorem.  Consider $\Pr \{ \ovl{\bs{\mcal{X}}}_n \leq \bs{z} \}$.  Let $t_\ell < 0$ for $\ell = 1, \cd,
k$. Note that \bel  \Pr \{ \ovl{\bs{\mcal{X}}}_n \leq \bs{z} \} & = & \Pr \li \{ \sum_{i=1}^n X_{\ell, i}
\leq n z_\ell, \; \ell = 1, \cd, k \ri \} \nonumber\\
&  = &  \Pr \li \{ \exp \li ( t_\ell \sum_{i=1}^n X_{\ell, i} \ri ) \geq \exp (n t_\ell z_\ell), \; \ell = 1, \cd, k \ri \} \nonumber\\
& \leq & \li [ \prod_{\ell = 1}^k \exp ( - n t_\ell z_\ell) \ri ] \; \bb{E} \li [  \prod_{\ell = 1}^k \exp \li ( t_\ell \sum_{i=1}^n X_{\ell, i}
\ri ) \ri ], \la{good88} \eel where (\ref{good88}) follows from multivariate Markov inequality.  Making use of (\ref{good88}) and Lemma
\ref{howgood}, we have \be \la{use88a}
  \Pr \{ \ovl{\bs{X}}_n \leq \bs{z} \}  \leq  \li [ \prod_{\ell = 1}^k \exp ( - n t_\ell z_\ell) \ri ] \;
\li [ \mu_0 + \sum_{\ell = 1}^k \mu_\ell \exp ( t_\ell ) \ri ]^n  =  \exp \li ( n \mscr{M} (t_1, \cd, t_k) \ri ) \ee for all $t_\ell \leq 0, \;
\ell = 1, \cd, k$, where $\mscr{M} (t_1, \cd, t_k)$ is defined by (\ref{def86636}).  Since $z_\ell \leq \mu_\ell$ for $\ell = 1, \cd, k$, we
have $z_0 \geq \mu_0$.  In view of the results of Lemma \ref{def68689},  we set \be  \la{use88b} t_\ell = \ln \f{z_\ell \mu_0}{z_0 \mu_\ell },
\qqu \ell = 1, \cd, k. \ee Then, \be \la{use88c} t_\ell  \leq 0, \qqu \ell = 1, \cd, k \ee and \bel \mscr{M} (t_1, \cd, t_k) & = & - \sum_{\ell
= 1}^k t_\ell z_\ell +  \ln \li [ \mu_0 + \sum_{\ell = 1}^k
\mu_\ell \exp ( t_\ell ) \ri ] \nonumber\\
& =  & - \sum_{\ell = 1}^k  z_\ell \ln \f{z_\ell  \mu_0 } { z_0 \mu_\ell } + \ln \li ( \sum_{\ell = 0}^k
\f{z_\ell \mu_0 } { z_0 } \ri ) \nonumber\\
& =  & - \sum_{\ell = 0}^k  z_\ell \ln \f{z_\ell  \mu_0 } { z_0 \mu_\ell } + \ln \li ( \f{ \mu_0 } { z_0 } \sum_{\ell = 0}^k  z_\ell \ri ) \nonumber\\
& =  & \sum_{\ell = 0}^k  z_\ell \ln \f { \mu_\ell  }{z_\ell  }.  \la{use88d} \eel  By virtue of (\ref{use88a}),  (\ref{use88b}), (\ref{use88c})
and (\ref{use88d}), we have \bee \Pr \{ \ovl{\bs{\mcal{X}}}_n \leq \bs{z} \}
& \leq &  \exp \li ( n \sum_{\ell = 0}^k  z_\ell \ln \f { \mu_\ell  }{z_\ell  } \ri ) \\
& = & \prod_{\ell = 0}^k \li ( \f { \mu_\ell  }{z_\ell  } \ri )^{n z_\ell}   \eee provided that $\bs{z} \leq \bs{\mu}$. In a similar manner, we
can show that \bee \Pr \{ \ovl{\bs{\mcal{X}}}_n \geq \bs{z} \} \leq \prod_{\ell = 0}^k \li ( \f { \mu_\ell  }{z_\ell  } \ri )^{n z_\ell} \eee
provided that $\bs{z} \geq \bs{\mu}$.  This completes the proof of the theorem.

\end{document}